\documentclass[12pt,reqno]{amsart}
\usepackage{mathrsfs}
\usepackage{amsfonts}

\usepackage{amssymb,amsmath, amsfonts, latexsym}
\usepackage{enumerate}
\newtheorem{thm}{Theorem}[section]
\newtheorem{cor}[thm]{Corollary}
\newtheorem{lem}[thm]{Lemma}
\newtheorem{prop}[thm]{Proposition}
\newtheorem{example}[thm]{Example}
\newtheorem{remarks}[thm]{Remark}
\newtheorem{defn}[thm]{Definition}

 \numberwithin{equation}{section}

\title[]{Centered Sobolev inequality and exponential convergence in $\Phi$-entropy}

\author{Lingyan Cheng}
\address{Lingyan Cheng. Institute of Applied Mathematics, Academy of Mathematics and Systems Science,
Chinese Academy of Sciences, 100190, Beijing China. }
\email{chengly@amss.ac.cn}

\author{Liming Wu}
\address{Liming Wu. Laboratoire de Math\'ematiques Appliqu\'ees, CNRS-UMR
6620, Universit\'e Blaise Pascal, 63177 Aubi\`ere, France. }
\email{Li-Ming.Wu@math.univ-bpclermont.fr}

\newcommand{\dd}{\mathbb{D}}

\newcommand{\rr}{\mathbb{R}}
\newcommand{\pp}{\mathbb{P}}

\def\AA{\mathcal A}
\def\BB{\mathcal B}

\def\DD{\mathcal D}
\def\FF{\mathcal F}

\def\EE{\mathcal E}
\def\HH{\mathcal H}

\def\LL{\mathcal L}

\def\vep{\varepsilon}
\def\<{\langle}
\def\>{\rangle}

\def\beq{\begin{equation}}
\def\neq{\end{equation}}

\def\bthm{\begin{thm}}
\def\nthm{\end{thm}}

\def\bprop{\begin{prop}}
\def\nprop{\end{prop}}

\def\brmk{\begin{remarks}}
\def\nrmk{\end{remarks}}

\def\bexa{\begin{example}}
\def\nexa{\end{example}}

\def\blem{\begin{lem}}
\def\nlem{\end{lem}}

\def\bcor{\begin{cor}}
\def\ncor{\end{cor}}

\def\bdef{\begin{defn}}
\def\ndef{\end{defn}}

\def\bexe{\begin{exe}}
\def\nexe{\end{exe}}

\def\bprf{\begin{proof}}
\def\nprf{\end{proof}}

\def\bdes{\begin{description}}
\def\ndes{\end{description}}

\def\benu{\begin{enumerate}}
\def\nenu{\end{enumerate}}

\begin{document}

\begin{abstract} In this short paper we find that the Sobolev
inequality $$\frac 1{p-2}\left[\left(\int f^{p} d\mu\right)^{2/p} - \int f^2 d\mu\right]
\le C \int |\nabla f|^2 d\mu$$ ($p\ge 0$) is equivalent to the
exponential convergence of the Markov diffusion semigroup $(P_t)$ to
the invariant measure $\mu$, in some $\Phi$-entropy. We provide the estimate of the exponential convergence in total variation and a bounded perturbation result under the Sobolev inequality. Finally in the one-dimensional case we get some two-sided estimates of the Sobolev constant by means of the generalized Hardy inequality.


\end{abstract}

\maketitle

\section{Introduction}
 \subsection{Centered Sobolev inequality} Let $\mu$ be a probability measure on some Polish space $E$
 equipped  with the Borel $\sigma$-field $\BB$. The main object of
 this paper is the following centered version of Sobolev inequality
\beq\label{Sobolev}
\frac 1{p-2}\left[\left(\int_E f^{p} d\mu\right)^{2/p} -
\int_E f^2 d\mu\right] \le C_S(p) \EE[f], \ 0\le f\in \dd(\EE)
\neq
where $p\in [0,+\infty)$,  $\EE$ is a conservative Dirichlet form on
$L^2(E,\mu)$ with domain $\dd(\EE)$ and $C_S(p)$ is the best
constant. This inequality will be denoted by $(S_p)$.

When $p=1$, this becomes the usual Poincar\'e inequality
\beq\label{Poincare}
{\rm Var}_\mu(f):=\mu(f^2) - [\mu(f)]^2\le
C_S(1) \EE[f], \ f\in \dd(\EE)
\neq
where $\mu(f):=\int_E f d\mu$.
Thus $C_S(1)$ is exactly the best Poincar\'e constant $C_P$.

When $p=2$, the left-hand side (LHS in short) of (\ref{Sobolev}),
understood as the limit when $p\to 2$, equals to $\frac 12 H(f^2)$,
where
$$
H(f)=\mu(f\log f) -\mu(f)\log \mu(f)
$$
is the entropy of $f$. So the Sobolev-type inequality
(\ref{Sobolev}) becomes
\beq\label{logS}
H(f^2)\le 2 C_S(2)\EE[f], \
0\le f\in\dd(\EE)
\neq
the usual log-Sobolev inequality (\cite{Ba94}). Thus
$C_S(2)$ coincides with the best log-Sobolev constant $C_{LS}$.

When $p>2$, $(S_p)$ is a centered version of the classic defective
Sobolev inequality :
\beq\label{defS} \left(\int_E f^p
d\mu\right)^{2/p}\le  A \EE[f] + B \int_E f^2 d\mu,\ 0\le f\in\dd(\EE).
\neq
For example
for the Lebesgue measure $\mu$ on $E=\rr^n$ and $\EE(f)=\int_E |\nabla
f|^2 d\mu$, the above Sobolev inequality holds with $B=0$, for  $p=
\frac{2n}{n-2}$ ($n>2$, see Aubin \cite{Au}). Notice that the defective Sobolev
inequality with $B=0$ fails for probability measure $\mu$.

The centered Sobolev inequality $(S_p)$ was studied by Aubin \cite{Au}
and Beckner \cite{Be} for the normalized volume measure on the unit
sphere $S^n$ in $\rr^{n+1}$. They obtained the exact result :
$C_S(p)=C_P=\frac 1n$ if $2<p\le \frac{2n}{n-2}$ (for $n\ge 3$).
Bakry and Ledoux \cite{BaLe}, using the diffusion semigroup method, proved the following sharp and general result of (see also Ledoux \cite[Theorem 3.1]{Led02}) :

\bthm\label{thm11} (\cite{BaLe})
Let $L$ be a Markov diffusion generator satisfying the Bakry-Emery's
curvature-dimension condition $CD(R,n)$ for some $R>0$ and
$n>2$. Then for every $1 \le p \le \frac{2n}{n-2}$, \eqref{Sobolev} holds with
$C_S(p) \le \frac{n-1}{nR}$.
\nthm

This deep theorem of Bakry-Ledoux generalizes the famous  Lichrowicz bound about $C_S(1)=C_P$.

When $p=0$, the LHS of (\ref{Sobolev}), understood as the limit when
$p\to0^+$, equals to $\frac 12[\mu(f^2)-e^{\mu(\log f^2)}]$. Setting
$f^2=e^g$, we see that $(S_0)$ becomes \beq\label{S0} \mu(e^g) -
e^{\mu(g)} \le 2C_S(0) \EE[e^{g/2}], \ g\in\dd(\EE)\bigcap
L^\infty(\mu). \neq

Relationship between the Sobolev inequalities for different $p$ is
summarized in

 \bthm\label{thmA} \benu[(a)]
\item For any $p\in \rr^+=[0,+\infty)$, $C_S(p)\ge C_S(1)=C_P$.
\item $pC_S(p)$ is nondecreasing in $p\in \rr^+$.
\item For any $p\in (0,2)$, the Sobolev inequality $(S_p)$ is
equivalent to the Poincar\'e inequality, more precisely
$$
\aligned &C_S(1)\le C_S(p) \le \frac{C_S(1)}{p},\ p\in(0,1);\\
&C_S(1)\le C_S(p) \le \frac{C_S(1)}{2-p},\ p\in(1,2).
\endaligned
$$
\nenu \nthm
This result is essentially contained in Bakry and Ledoux \cite{BaLe}.

 In other words this family of Sobolev inequalities for
different $p$ has four interesting cases : (1) $p=0$ ; (2) $p=1$ ;
(3) $p=2$ and (4) $p>2$.

\subsection{Semigroup}
Let $(P_t)$ be a Markov semigroup such that $\mu P_t=\mu$ for all $t\ge0$
(i.e. $\mu$ is an invariant measure), strongly continuous on $L^2(\mu)$. Let $\LL$ be the generator of
$(P_t)$, whose domain in $L^p(\mu):=L^p(E,\BB,\mu)$ is denoted by
$\dd_p(\LL)$ ($1\le p<\infty$). We always assume that

 {\it {\bf (A1)} $\dd_2(\LL)$ is contained in $\dd(\EE)$ and dense in
 $\dd(\EE)$ w.r.t. the norm $\sqrt{\mu(f^2) +\EE[f]}$ (i.e. $\dd_2(\LL)$ is a form core of
 $\EE$), and
 $$
\int f(-\LL f) d\mu = \EE[f],\ f\in\dd_2(\LL).
 $$}
In other words $\EE$ is the symmetrized Dirichlet form of $\LL$. This assumption holds automatically if $\LL$ is self-adjoint (i.e. $(P_t)$ is symmetric on $L^2(\mu)$).

It is well known that the Poincar\'e inequality $(S_1)$ is
equivalent to the exponential convergence of $P_t$ to $\mu$ in
$L^2(\mu)$ :
$$
{\rm Var}_\mu[P_tf] \le e^{-2t/C_S(1)} {\rm Var}_\mu[f],\ t>0, \ f\in
L^2(\mu).
$$
And if $(P_t)$ is a diffusion semigroup, the log-Sobolev inequality
$(S_2)$ is equivalent to the exponential convergence of $P_t$ to
$\mu$ in the relative entropy
$$
H(P_tf)\le e^{-2t/C_{S}(2)}H(g), \ t>0,\ 0\le g\in L^1(\mu).
$$
See Bakry \cite{Ba94}. Notice that the later equivalence is false in
the jump case (\cite{wu00h}).

But unlike Poincar\'e and log-Sobolev, the role of the Sobolev
inequality (\ref{Sobolev}) for $p$ different from $1,2$ in the
exponential convergence of $P_t$ is unknown. Our first purpose of
this paper is to fill this gap.

This paper is organized as follows.
In the next section we establish the equivalence between the Sobolev
inequality and the exponential convergence of $P_t$ to $\mu$, in
some $\Phi$-entropy sense. Several corollaries and applications are
derived for illustrating the usefulness of the Sobolev inequality
(\ref{Sobolev}), especially for the rate of the exponential
convergence of $P_t$ to $\mu$ in total variation.

In \S3 we recall the relationship between the defective Sobolev inequality and centered Sobolev inequality when $p>2$ and present a bounded perturbation result.

In \S4 we present some two-sided estimates of the optimal constant $C_S(p)$ of Sobolev inequality when $p>2 $ on the real line, by following Barthe and Roberto \cite{BR03}.

\section{Equivalence between Sobolev inequality and exponential convergence}

\subsection{\bf Framework} Besides {\bf (A1)}, we assume

{\it {\bf (A2) (Existence of the carr\'e-du-champs operator)} there
is an algebra $\AA$ contained in $\dd_2(\LL)$ and dense in
$\dd(\EE)$ w.r.t. the norm $\|f\|_{2,1}:=\sqrt{\mu(f^2)+\EE[f]}$. So
the {\it carr\'e-du-champs} operator
$$
\Gamma(f,g):= \frac 12[\LL(fg)-f\LL g-g\LL f], f,g\in \AA
$$
is well defined. $\Gamma(f,g)$ can be extended as a continuous
mapping from $\dd(\EE)\times \dd(\EE)\to L^1(\mu)$.

{\bf (A3)} $(P_t)$ is a diffusion semigroup, i.e. $P_t$ is the
transition probability semigroup of a {\bf continuous} Markov
process $(X_t)$ valued in $E$ defined on $(\Omega, \FF_t, \pp_\mu)$.
}

Under those assumptions, for every $f\in \dd_2(\LL)$,
$$
M_t(f):=f(X_t)-f(X_0) - \int_0^t \LL f(X_s) ds
$$
is a $L^2(\pp_\mu)$-martingale, and
$$
\<M(f), M(g)\>_t=2\int_0^t \Gamma(f,g) (X_s) ds
$$
(this holds at first for $f,g\in \AA$, then for $f\in\dd_2(\LL)$ by
continuous extension). Consequently if $f_1,\cdots, f_n\in
\dd_\infty(\LL)=\{f\in \dd_2(\LL)\ ;\ f,\ \LL f\in L^\infty(\mu)\}$ and
$F:\rr^n\to \rr$ is infinitely differentiable, then by Ito's formula
$F(f_1,\cdots,f_n)\in \dd_2(\LL)$ and
\beq\label{11}
\LL F(f_1,\cdots, f_n)= \sum_{i=1}^n
\partial_i F(f_1,\cdots,f_n)\LL f_i +\sum_{i,j=1}^n
\partial_i\partial_j F(f_1,\cdots,f_n)\Gamma(f_i,f_j).
\neq
We write $\Gamma[f]=\Gamma(f,f)$. For any $C^\infty$-function
$\Phi$ on $\rr$, integrating (\ref{11}) we have
\beq\label{21}
\int \Phi'(f) (-\LL f) d\mu= \int \Phi^{\prime\prime}(f)\Gamma[f]d\mu, \
 f\in\dd_\infty(\LL).
\neq
Next (\ref{11}) implies that $\Gamma$ is
a derivation,
\beq\label{22}
\Gamma(\Phi(f),g)  = \Phi^{\prime}(f)\Gamma(f,g), \ f\in\dd_\infty(\LL), \ g\in
\dd_2(\LL).
\neq

\subsection{Exponential convergence in the $\Phi$-entropy}

\begin{defn}\label{def-entropy} Given a lower bounded convex function $\Phi: \rr\to
(-\infty,+\infty]$, the $\Phi$-entropy of a function $f\in L^1(\mu)$
is defined as
$$
H_\Phi^\mu(f)=\mu(\Phi(f)) - \Phi(\mu(f)).
$$
\end{defn}

The main result of this section is

\bthm\label{thm1}
For the diffusion Markov semigroup $(P_t)$ with
invariant probability measure $\mu$ satisfying {\bf (A1), (A2) and
(A3)}, the Sobolev inequality (\ref{Sobolev}) is equivalent to the
exponential convergence in the $\Phi$-entropy
\beq\label{thm1a}
H_\Phi^\mu(P_tf)\le e^{-\frac{2t}{C_S(p)}}H_\Phi^\mu(f), \ f\in L^1(\mu)
\neq
where
\beq\label{thm1b}
\Phi(x)=\begin{cases} -|x|^{2/p},\quad &\text{ if }\ p\in(2,+\infty);\\
|x|\log |x|, &\text{ if }\ p=2;\\
|x|^{2/p}, &\text{ if }\ p\in(0,2);\\
e^x, &\text{ if }\ p=0.
\end{cases}
\neq
\nthm We begin with a known result (see Chafai \cite{Chafai}). \blem\label{lem23}
Let $\Phi$ be a lower bounded $C^2$-convex function and $\DD$ be a
class of functions in $\dd_\infty(\LL)$, stable for $(P_t)$ (i.e. if
$f\in\DD$, $P_tf\in \DD$). The exponential convergence in the
$\Phi$-entropy
$$ H_\Phi^\mu(P_tf)\le e^{-\frac{2t}{C(\Phi)}}H_\Phi^\mu(f), \ f\in
\DD$$ for some positive constant $C(\Phi )$ is equivalent to
\beq\label{Phi<} H_\Phi^\mu(f) \le \frac{C(\Phi)}{2} \int
\Phi^{\prime\prime}(f)\Gamma[f] d\mu,\ f\in\DD.\neq \nlem

\bprf Since for $f\in\DD$
$$\aligned
\frac d{dt} H_\Phi^\mu(P_tf)& = \int \Phi'(P_tf) \LL P_tf d\mu\\
&= - \int \Phi^{\prime\prime}(P_tf) \Gamma[P_tf] d\mu
\endaligned
$$
by (\ref{21}), the equivalence above follows from Gronwall's lemma.
\nprf

\bprf[\bf Proof of Theorem \ref{thm1}] For the exponential
convergence in the $\Phi$-entropy we may restrict to $f\in
\DD=\{f\in \dd_\infty(\LL); \ \exists \vep>0, f\ge \vep \}$. In that
case as $\Phi$ is $C^2$ on $(0,+\infty)$, we can apply Lemma
\ref{lem23}.

At first this equivalence is well known for $p=1,2$ as recalled in
the Introduction.

We begin with the case $p>2$. By Lemma \ref{lem23},
 the exponential convergence (\ref{thm1a}) is equivalent to
$$
[\mu(f)]^{2/p} - \mu(f^{2/p})\le C_S(p)\frac{p-2 }{p^2} \int_E
f^{\frac{2-2p}{p}}\Gamma[f] d\mu,\ f\in\DD.
$$
Setting $h=f^{1/p}$, this last inequality is equivalent to
$$
[\mu(h^p)]^{2/p} - \mu(h^2)\le C_S(p)(p-2) \int_E \Gamma[h]
d\mu,\ \vep^{1/p}\le h\in\dd_\infty(\LL)
$$
which is exactly the Sobolev inequality (\ref{Sobolev}).

For $p\in (0,2)$, by Lemma \ref{lem23},
 the exponential convergence (\ref{thm1a}) is equivalent to
$$
 \mu(f^{2/p})-[\mu(f)]^{2/p} \le C_S(p)\frac{2-p}{p^2} \int_E
f^{\frac{2-2p}{p}}\Gamma[f] d\mu,\ f\in\DD.
$$
Setting $h=f^{1/p}$, this last inequality is equivalent to
$$
\mu(h^2)-[\mu(h^p)]^{2/p}  \le C_S(p)(2-p) \int_E \Gamma[h] d\mu,\
\vep^{1/p}\le h\in\dd_\infty(\LL)
$$
which is exactly the Sobolev inequality (\ref{Sobolev}).

Finally for $p=0$, by Lemma \ref{lem23},
 the exponential convergence (\ref{thm1a}) is equivalent to
$$
 \mu(e^f)-e^{\mu(f)} \le \frac{C_S(0)}{2} \int_E
e^f\Gamma[f] d\mu,\ f\in\DD
$$
which is exactly the Sobolev inequality (\ref{S0}) for $p=0$.
 \nprf

\subsection{Exponential convergence in Hellinger metric} Now we present an application to the exponential
convergence in the Hellinger metric $d_\HH$. Recall that for two
probability measures $\nu=g d\alpha,\mu=fd\alpha$ where $\alpha$ is
some reference measure,
$$
d_\HH^2(\nu,\mu):=\int(\sqrt{g}-\sqrt{f})^2 d\alpha.
$$
$d_\HH$ is in fact independent of the choice of $\alpha$. \bcor
Assume that the adjoint operator $\LL^*$ of $\LL   $ satisfies also {\bf
(A1), (A2), (A3)}. The Sobolev inequality (\ref{Sobolev}) for $p=4$
is equivalent to
$$
d_\HH(P_t^* f\mu,\mu) \le e^{-t/C_S(4)} d_\HH(f\mu,\mu), \ t>0
$$
for any $\mu$-probability density function $f$. \ncor

Recall that the distribution of $X_t$ is $P_t^* f \mu$ if the initial distribution
of $X_0$ is $f\mu$.

\bprf We have for any $\mu$-probability density function $f$,
$$
d_\HH^2 (f\mu,\mu)=\int (\sqrt{f}-1)^2d\mu=2(1-\mu(\sqrt{f})).
$$
And for the exponential convergence in (\ref{thm1a}) (with $p=4$), one may restrict to the functions $f\ge0$ such that $\mu(f)=1$ by homogeneity. So this corollary follows directly by Theorem \ref{thm1}(a).
\nprf

\brmk
{\rm Let $\|\nu-\mu\|_{TV}:=\sup_{|f|\le 1}|\nu(f)-\mu(f)|$  (the total variation). It is known that (Gibbs and Su \cite{GS02})
$$
d_\HH^2(\nu,\mu) \le \|\nu-\mu\|_{TV}\ \text{ and }\
\|\nu-\mu\|_{TV} \le 2 d_\HH(\nu,\mu).
$$
So under the Sobolev inequality (\ref{Sobolev}) with $p=4$, we have
$$
\|P_t^*f \mu-\mu\|_{TV} \le 2 e^{-t/C_S(4)} d_\HH(f\mu,\mu)\le 2 \sqrt{2}
e^{-t/C_S(4)}
$$
an explicit estimate of the exponential convergence in total variation.
}
\nrmk

\subsection{Exponential convergence in total variation}
We now generalize the result above to general $p>2$ different from $4$.

\bcor Assume that $\LL^*$ satisfies {\bf (A1), (A2), (A3)}. If the
Sobolev inequality holds for some $p>2$, then for any $\mu$-probability
density $f$,
$$
\|P_t^* f \mu -\mu\|_{TV}\le 2p\sqrt{\frac{1}{p-2}}e^{-t/C_S(p)}
(1-\mu[f^{2/p}])^{1/2} \le 2p\sqrt{\frac{1}{p-2}} e^{-t/C_S(p)}.
$$
\ncor

\bprf It follows from Theorem \ref{thm1} and the lemma below.\nprf

\blem\label{lem21} Let $a\in (0,1)$. Then for any $f\ge0$ such that
$\mu(f)=1$, \beq\label{lem21a} 1-\mu(f^a)\le \frac{1}{2} \int
|f-1|d\mu \neq and \beq\label{lem21b} 1-\mu(f^a)\ge \frac{a(1-a)}{8}
\left(\int |f-1|d\mu\right)^2. \neq \nlem

\bprf We have
$$
\aligned 1-\mu(f^a)&= \int_{E} (1-f^a) d\mu \le  \int_{\{f<1\}}
(1-f^a)
d\mu\\
&\le\int_{\{f<1\}} (1-f)d\mu=\frac 12 \int_E |f-1| d\mu,\\
\endaligned
$$
that is (\ref{lem21a}).

For (\ref{lem21b}) we may assume that $\mu(f=1)<1$. Letting
$A=\{f<1\}$,
$$\bar f= \frac{\mu(f1_A)}{\mu(A)} 1_A+ \frac{\mu(f1_{A^c})}{\mu(A^c)} 1_{A^c}
$$
(which is the conditional expectation of $f$ knowing $\sigma(A)$), by Jensen's inequality we have
$$
1-\mu(f^a)\ge 1-\mu(\bar f^a),\ \int_E |f-1|d\mu =
2\int_{\{f<1\}}(1-f) d\mu=\int_E |\bar f-1|d\mu.
$$
So it is enough to prove (\ref{lem21b}) for $f=\bar f$, a two-valued
function. Let $x< y$ be the two values of $f$ (so $0\le x<1<y$), and
$$
\alpha:=\mu(f=x)=1-\mu(f=y)=:1-\beta.
$$
Since $\mu(f)=\alpha x+ \beta y=1$, $y=\frac{1-\alpha x}{\beta}$.
Consider
$$
h(x)=1-\mu(f^a)=1-[\alpha x^a+ \beta y^a].
$$
We have $h(1)=h'(1)=0$, and $\int |f-1|d\mu=2\alpha(1-x)$. Hence for
(\ref{lem21b}), by Taylor's formula we have only to show that
\beq\label{lem21c} \min_{x\in(0,1)}h^{\prime\prime}(x)\ge
a(1-a)\alpha^2. \neq Notice that
$$
\aligned
h^{\prime\prime}(x) &=-a(a-1)\alpha[x^{a-2}+\frac \alpha\beta y^{a-2}],\\
h^{\prime\prime\prime}(x) &=-a(a-1)(a-2)\alpha[x^{a-3}-(\frac \alpha\beta)^2 y^{a-3}]
\endaligned$$
and $h^{(4)}(x)>0$ for all $x\in(0,1)$. We now divide our discussion
into two cases.

{\bf Case 1. $\alpha\le 1/2$.} In this case,
$h^{\prime\prime\prime}(1)\le 0$, then $h^{\prime\prime\prime}(x)<0$
for all $x\in (0,1)$, consequently
$$h^{\prime\prime}(x)\ge
h^{\prime\prime}(1)=a(1-a)\alpha^2\frac 1{\alpha\beta}\ge 4a(1-a)
\alpha^2$$ which implies (\ref{lem21c}).\\
 {\bf Case 2.
$\alpha>1/2$.} Since $lim_{x\to0_+}h^{\prime\prime\prime}(x)=-\infty$
and $h^{\prime\prime\prime}(1)>0$, there is a unique $x_0\in(0,1)$
such that $h^{\prime\prime\prime}(x_0)=0$, i.e.
$x_0^{a-3}=(\alpha/\beta)^2y_0^{a-3}$ ($y_0=(1-\alpha x_0)/\beta$)
or $x_0=\frac 1\alpha
\frac{\left(\frac\alpha\beta\right)^{\frac{a-1}{a-3}}}{1+\left(\frac\alpha\beta\right)^{\frac{a-1}{a-3}}}$.
Consequently
$$
\min_{x\in(0,1]}
h^{\prime\prime}(x)=h^{\prime\prime}(x_0)=a(1-a)\alpha^2\left(\frac
1\alpha + \frac1\beta
\left(\frac\beta\alpha\right)^{\frac{2(a-2)}{a-3}}\right)x_0^{a-2}\ge
a(1-a)\alpha^2
$$
for $x_0^{a-2}>1$. The last bound is optimal because it becomes
equality if $\alpha\to1$. That completes the proof of
(\ref{lem21c}). \nprf

\section{Defective Sobolev inequality implies centered Sobolev inequality and a bounded perturbation result}
\subsection{Defective Sobolev inequality implies Sobolev inequality }

\bthm\label{thm31} (\cite{BGL14})
 If  the defective Sobolev inequality \eqref{defS} holds with some positive constants $A,B$ for some $p>2$,  and the Poincar\'{e} inequality \eqref{Poincare} holds with the best constant $C_P>0$,
then we have
\beq
\left(\int_E|f|^p d \mu\right)^{2/p}-\int_E f^2 d\mu \le  \left((p-1) A+C_P[(p-1)B-1]^+\right)\EE[f].
\neq
\nthm

The above theorem \ref{thm31} is a direct consequence of the following lemma.

\blem \label{lem43} Let $p>2$ and $f:E \to \rr$ be a square integrable function on a probability space $(E,\mu)$. Then for all $a \in \rr$, we have
\beq\label{a1}
\left(\int_E|f|^p d \mu\right)^{2/p}-\int_E f^2 d\mu \le (p-1)\left( \int_E |f-a|^p d\mu \right)^{2/p}-\int_E (f-a)^2 d\mu.
\neq
\nlem
This lemma is also contained in \cite{BGL14} and will be used in the next section.

Notice that if the defective Sobolev inequality holds for some $p>2$, then  $P_t(x,dy)=p_t(x,y)\mu(dy)$ with the density $p_t(x,y)$ bounded (\cite{BaLe}). That implies $P_t$ is a Hilbert-Schmidt operator, then compact on $L^2(\mu)$ : in particular the Poincar\'e inequality holds true.

\subsection{Bounded perturbation}

It is well known that $\Phi$-entropy  $H_\Phi^\mu(f)$ defined in definition \ref{def-entropy} has the following variational form :
\beq\label{variational form}
H_\Phi^\mu(f)=\mu(\Phi(f)) - \Phi(\mu(f))=\inf_{c \in \rr} \int_E \Phi(f) - \Phi(c) -\Phi^\prime(c) (f-c) d\mu
\neq
for all $f \in L^1(\mu)$. The following proposition shows that the Sobolev inequality \eqref{Sobolev} is stable by bounded transformation of the probability measure $\mu$.

\bprop\label{prop33} Assume that the Dirichlet form $\EE[f]=\int \Gamma[f]d\mu$ for some {\bf carr\'e-du-champs} operator $\Gamma$ which is a derivation, i.e. $\Gamma(\Phi(f), g)=\Phi'(f) \Gamma(f,g)$ for all $f,g\in\dd(\EE)\bigcap L^\infty(\mu)$ and $\Phi\in C^1(\rr)$. Assume that the probability measure $\mu$ satisfies Sobolev inequality \eqref{Sobolev} with the best constant $C_S(p)$ for $p\ge 0$. Let $\tilde{\mu}$ be the probability measure defined by $d \tilde{\mu} = \frac{1}{Z}e^{-V(x)}d \mu$ such that ${\bf Osc}(V):=\sup_{x,y \in E}| V(x)-V(y)|< +\infty$, where $Z>0$ is the normalization constant. Then $ \tilde{\mu}$ satisfies Sobolev inequality
$$
\frac 1{p-2}\left[\left(\int_E f^{p} d\tilde \mu\right)^{2/p} -
\int_E f^2 d\tilde\mu\right] \le  e^{{\bf Osc}(V)} C_S(p) \int \Gamma[f]d\tilde \mu, \ 0\le f\in \dd(\EE) \bigcap L^\infty(\mu).
$$
\nprop

\bprf
By the proof of Theorem \ref{thm1}, the Sobolev inequality \eqref{Sobolev} is equivalent to
$$
H_\Phi^\mu(f)\le \frac{C_S(p)}{2}\int_E \Phi^{\prime \prime} (f) \Gamma(f) d\mu, \ f\in\DD:=\{g\in \dd(\EE); \exists \vep>0; \vep\le g\le 1/\vep\}.
$$
where
$$
\Phi(x)= \begin{cases}-|x|^{2/p}, &\mbox{if} \ p \in (2,+\infty);\\
|x| \log|x|, &\mbox{if} \ p =2;\\
|x|^{2/p}, &\mbox{if} \ p \in (0,2);\\
e^x,   &\mbox{if} \ p =0.
\end{cases}
$$
We have by \eqref{variational form},
$$
\begin{aligned}
H_\Phi^{\tilde{\mu}}(f)&=\inf_{c \in \rr} \int_E \Phi(f) - \Phi(c) -\Phi^\prime(c) (f-c) d \tilde{\mu}\\
&=\inf_{c \in \rr} \int_E [\Phi(f) - \Phi(c) -\Phi^\prime(c) (f-c)] \frac{1}{Z}e^{-V} d \mu\\
&\le \frac{1}{Z}\exp\left(-\displaystyle\inf_{x\in E} V(x)\right) H_\Phi^{\mu}(f)\\
&\le \frac{1 }{2Z}\exp\left(-\displaystyle\inf_{x\in E} V(x)\right)C_S(p)\int_E \Phi^{\prime \prime} (f) \Gamma(f)d \mu\\
&=\frac{ 1}{2Z}\exp\left(-\displaystyle\inf_{x\in E} V(x)\right) C_S(p)\int_E \Phi^{\prime \prime} (f) \Gamma(f) Z e^{V} d \tilde{\mu}\\
&\le \frac{1}{2}\exp\left(\displaystyle\sup_{x\in E} V(x)-\displaystyle\inf_{x\in E} V(x)\right)C_S(p)\int_E \Phi^{\prime \prime} (f) \Gamma(f) d \tilde{\mu}\\
&= \frac{1}{2}e^{{\bf Osc}(V)} C_S(p) \int_E \Phi^{\prime \prime} (f) \Gamma(f) d \tilde{\mu}
\end{aligned}
$$
which implies the result.
\nprf

\subsection{Reflected Brownian motion}
Given a domain $\Omega$ of $\rr^d$,
let $W^{1,p}(\Omega) $ be the Sobolev space of the functions on $\Omega$ with the norm $ \| f\|_{ W^{1,p}(\Omega)}= \left(\int_\Omega (|\nabla f|^p +|f|^p) dx\right)^{\frac{1}{p}}$. Recall the extension theorem on Sobolev space:
\bthm\label{extension thm}
Suppose that $\Omega \subset \rr^d$ is a bounded domain with Lipschitz boundary. Then there exist a bounded linear operator $L: W^{1,p}(\Omega) \ni u \to v \in W^{1,p}(\rr^d)$ and a constant $C>0$ such that
\benu
\item
$ v(x)=u(x)$ for a.e. $x \in \Omega$;
\item
$\|v\|_{W^{1,p}(\rr^d)} \le C\|u\|_{ W^{1,p}(\Omega)}$.
\nenu
\nthm
According to the well known Sobolev inequality on $\rr^d$, we have the following
\bcor\label{cor36}
For any bounded domain $ \Omega \in \rr^d$ ($d\ge 2$) with Lipschitz boundary, the Sobolev inequality \eqref{Sobolev} holds for $ u \in W^{1,2}(\Omega)$ with the normalized Lebesgue measure $\mu(dx)=\frac{dx}{Vol(\Omega)}$ on $\Omega$ for any $ p\in (2,  \frac{2d}{d-2}]$ (this last quantity is interpreted as $+\infty$ if $d=2$).
\ncor
\bprf
By Theorem \ref{extension thm}, we have for all $ u\in W^{1,2}(\Omega)$,
$$
\begin{aligned}
\left(\int_\Omega |u|^p dx \right)^{\frac{2}{p}} &\le \left(\int_{\rr^d} |v|^p dx\right)^{\frac{2}{p}}\\
&\le C(d,p) \int_{\rr^d} |\nabla v|^2 dx\\
&\le C(d,p) C \int_\Omega (|\nabla u|^2 + |u|^2) dx,
\end{aligned}
$$
where $v= L u$,  $C(d,p)$ is the best Sobolev constant. Then
the defective Sobolev inequality \eqref{defS} holds with $A=B=C(d,p) C$. The result follows by Theorem \ref{thm31}.
\nprf

\section{Sobolev inequality in dimension one}
In \cite{BR03}, F.~Barthe and C.~Roberto provide the estimate of the optimal constant of Sobolev inequality when $1< p \le 2$ on the real line. In this section we generalize the estimate of the optimal constant to the case $p>2$ on the real line by the method in \cite{BR03}.

\bthm \label{estimate}
Let $p>2$ and $\mu,\ \nu$ (non-negative) be Borel measures on $\rr$ with $\mu (\rr)=1$ and $d\nu(x)=n(x)dx$, where $n(x)dx$ is the absolutely continuous component of $\nu$. Let $m$ be a median of $\mu$. Let $C>0$ be the optimal constant satisfying :
\beq \label{41}
\left(\int_\rr|f|^p d \mu\right)^{2/p}-\int_\rr f^2 d\mu \le C \int_\rr f'^2d \nu
\neq
for every smooth function $f:\ \rr \to \rr$.

Then we have $\max (b_-(p),b_+(p)) \le C \le 4 \max (B_-(p),B_+(p))$, where
\begin{align*}
&b_+(p)=\sup_{x>m} \{ \mu([x,+\infty))\bigg[\bigg(1+\frac{1}{2\mu[x,+\infty)}\bigg)^{(p-2)/p}-1\bigg] \int_m^x \frac{1}{n(t)}dt \};\\
&b_-(p)=\sup_{x<m} \{ \mu((-\infty,x])\bigg[\bigg(1+\frac{1}{2\mu(-\infty,x]}\bigg)^{(p-2)/p}-1\bigg] \int_x^m \frac{1}{n(t)}dt \};\\
&B_+(p)=\sup_{x>m} \{ \mu([x,+\infty))\bigg[\bigg(1+\frac{(p-1)^{p/(p-2)}}{\mu[x,+\infty)}\bigg)^{(p-2)/p}-1\bigg] \int_m^x \frac{1}{n(t)}dt \};\\
&B_-(p)=\sup_{x<m} \{ \mu((-\infty,x])\bigg[\bigg(1+\frac{(p-1)^{p/(p-2)}}{\mu(-\infty,x]}\bigg)^{(p-2)/p}-1\bigg] \int_x^m \frac{1}{n(t)} dt \}.
\end{align*}
\nthm

We will use the following Proposition and Lemmas to prove Theorem \ref{estimate}.

\bprop\label{prop42}{\em (See \cite{BR03})}
Let $\mu,\ \nu$ (non-negative) be Borel measures on $[m, \infty)$, where $m$ is a median of $\mu$ and $d\nu(x)=n(x)dx$, where $n(x)dx$ is the absolutely continuous component of $\nu$. Let $G$ be a family of non-negative Borel measurable functions on $[m,\infty)$. We set $\phi (f)=\sup_{g\in G} \int_m^\infty f g d\mu $ for any measurable function $f$. Let $A$ be the smallest constant such that for every smooth function $f$ with $f(m)=0$, we have
$$
\phi(f^2) \le A \int_m^\infty f'^2d \nu.
$$
Then $B \le A \le 4B$, where
$$
B=\sup_{x>m}\phi(1_{[x,\infty)})\int_m^x \frac{dt}{n(t)}.
$$
\nprop

\blem\label{lem44}
Let $\varphi$ be a non-negative integrable function on a probability space $(E,\mu)$. Let $A > 0$ and $a>1$ be some constants, then we have
$$
\aligned
&\quad A\left( \int \varphi^a d\mu \right)^{1/a}-\int \varphi d\mu \\
&=\sup \left\{ \int \varphi g d\mu\ ;\  g\ge-1 \ and \ \int (g+1)^{a/(a-1)} d\mu \le A^{a/(a-1)}\right\}\\
&\le \sup \left\{ \int \varphi g d\mu\ ;\  g \ge 0 \ and\ \int (g+1)^{a/(a-1)} d \mu \le A^{a/(a-1)}+1\right\}.
\endaligned
$$
\nlem

\bprf
For any Borel measurable function $h\ge0$, by H\"{o}lder's inequality, we have
$$
\left( \int \varphi^a d\mu \right)^{1/a}= \sup \left\{ \int \varphi h d\mu\ ;\ h\ge0\ and\ \int h^{a/(a-1)} d\mu \le 1 \right\}.
$$
Hence
\beq \label{43}
A \left( \int \varphi^a d\mu \right)^{1/a}= \sup \left\{ \int \varphi h d\mu\ ;\ h\ge 0\ and\ \int h^{a/(a-1)} d\mu \le A^{a/(a-1)} \right\}.
\neq
Using (\ref{43}), we have
$$
\aligned
&\quad A \left( \int \varphi^a d\mu \right)^{1/a} - \int \varphi d\mu\\
&=\sup\left\{ \int \varphi h d\mu-\int \varphi d\mu\ ;\ h\ge 0\ and \ \int h^{a/(a-1)} d\mu \le A^{a/(a-1)} \right\}\\
&=\sup\left\{ \int \varphi (h-1) d\mu\ ;\ h\ge 0\ and \ \int h^{a/(a-1)} d\mu \le A^{a/(a-1)} \right\}\\
&=\sup\left\{ \int \varphi g d\mu \ ;\ g\ge -1\ and \ \int (g+1)^{a/(a-1)} d\mu \le A^{a/(a-1)} \right\}\\
&\le \sup\left\{ \int \varphi g 1_{g \ge 0}d\mu \ ;\ g\ge-1\ and \ \int (g+1)^{a/(a-1)} d\mu \le A^{a/(a-1)} \right\}\\
&\le \sup\left\{ \int \varphi g d\mu \ ;\ g \ge 0 \ and \ \int (g+1)^{a/(a-1)} d\mu \le A^{a/(a-1)}+1 \right\}.
\endaligned
$$
The last inequality is derived by
$$
\int (g 1_{g \ge 1}+1)^{a/(a-1)} d\mu =\int (g+1)^{a/(a-1)}1_{g \ge 0}d\mu+\mu(g<0)\le A^{a/(a-1)}+1.
$$
Hence the lemma is established.
\nprf

\blem \label{lem45}
Let $a > 1$, $\mu$ be a finite measure on $X$. Let $A \subset X$ be a measurable subset with $\mu(A) >0$ and K be a constant with $K> \mu(X)$. Then we have
$$
\aligned
&\quad \sup \left\{ \int_X 1_A g d\mu \ ;\ g\ge 0\ and \ \int_X (g+1)^{a/(a-1)} d\mu \le K \right\}\\
&=\mu(A)\left[ \left( 1+\frac{K-\mu(X)}{\mu(A)}\right)^{(a-1)/a}-1\right].
\endaligned
$$
\nlem

\bprf
Simply, we denote by $S$ the right hand side of the above equality. Without loss of generality, we can assume $g=0$ on $A^c$, hence
$$
S= \sup \left\{ \int_A g d\mu \ ;\ g \ge 0\ and\ \int_A (g+1)^{a/(a-1)} d\mu +\mu(A^c) \le K \right\}.
$$
For any $g \ge 0\ and\ \int_A (g+1)^{a/(a-1)} d\mu +\mu(A^c) \le K $, by Jensen's inequality, we have
$$
\aligned
\left( 1+\int_A g \frac{d\mu}{\mu (A)}\right)^{a/(a-1)}&=\left( \int_A (1+g) \frac{d\mu}{\mu (A)} \right)^{a/(a-1)} \le \int_A (1+g)^{a/(a-1)} \frac{d\mu}{\mu (A)}\\
&\le \frac{K-\mu(A^c)}{\mu(A)}=1+\frac{K-\mu(X)}{\mu(A)}.
\endaligned
$$
Hence
\beq \label{44}
\int_A g d \mu \le \mu(A) \left[ \left( 1+\frac{K-\mu(X)}{\mu(A)}\right)^{(a-1)/a}-1\right].
\neq
We take $g= \left( 1+\frac{K-\mu(X)}{\mu(A)}\right)^{(a-1)/a}-1$, then equality in (\ref{44}) holds. Hence
$$
S=\mu(A) \left[ \left( 1+\frac{K-\mu(X)}{\mu(A)}\right)^{(a-1)/a}-1\right]
$$
which is the desired result.
\nprf

Now we prove Theorem \ref{estimate}.

\bprf[\bf Proof of Theorem \ref{estimate}.] {\bf Step 1.} We estimate the upper bound of $C$. For any smooth function $f:\rr \to \rr$, let $F=f-f(m)$, $F_+=F1_{(m,\infty)} $ and $F_-=F1_{(-\infty,m)}$. It is easy to see they are all continuous and $F^2=F^2_++F^2_-$, $|F|^p=|F_+|^p+|F_-|^p$ when $p>2$. We set $A=p-1>0$ and $a=\frac{p}{2}>1$, then $\frac{a}{a-1}=\frac{p}{p-2}$. By Lemma \ref{lem43} and Lemma \ref{lem44}, we have
$$
\aligned
&\quad \left( \int |f|^p d\mu  \right)^{2/p}-\int f^2 d\mu\\
&\le (p-1)\left( \int |f-f(m)|^p d\mu \right)^{2/p}-\int (f-f(m))^2 d\mu\\
&= (p-1)\left( \int |F|^p d\mu \right)^{2/p}-\int F^2 d\mu\\
&=(p-1)\left( \int |F_+|^p d\mu \right)^{2/p}-\int F_+^2 d\mu+(p-1)\left( \int |F_-|^p d\mu \right)^{2/p}-\int F_-^2 d\mu\\
&\le \sup \left\{ \int F_+^2 g d\mu \ ;\  g \ge 0 \ and\ \int (g+1)^{p/(p-2)} d \mu \le (p-1)^{p/(p-2)}+1\right\}\\
&\quad +\sup \left\{ \int F_-^2 g d\mu \ ;\  g \ge 0 \ and\ \int (g+1)^{p/(p-2)} d \mu \le (p-1)^{p/(p-2)}+1\right\}.
\endaligned
$$
Now we deal with $F_+$. Since $F_+=0$ on $(-\infty,m]$, we have by Proposition \ref{prop42},
$$
\aligned
&\quad (p-1)\left( \int |F_+|^p d\mu \right)^{2/p}-\int F_+^2 d\mu\\
&\le \sup \left\{ \int F_+^2 g d\mu \ ;\  g \ge 0 \ and\ \int (g+1)^{p/(p-2)} d \mu \le (p-1)^{p/(p-2)}+1\right\}\\
&\le 4 \tilde{B}_+(p) \int F_+'^2 d\nu
\endaligned
$$
where
$$
\aligned
&\quad \tilde{B}_+(p)= \\
&\sup_{x>m} \left[\sup \left\{ \int 1_{[x,\infty)} g d\mu \ ;\  g \ge 0 \ and\ \int (g+1)^{p/(p-2)} d \mu \le (p-1)^{p/(p-2)}+1\right\}\int_m^x \frac{1}{n(t)}dt\right].
\endaligned
$$
By Lemma \ref{lem45}, we have
$$
\tilde{B}_+(p)=\sup_{x>m} \{ \mu([x,+\infty))\bigg[\bigg(1+\frac{(p-1)^{p/(p-2)}}{\mu[x,+\infty)}\bigg)^{(p-2)/p}-1\bigg] \int_m^x \frac{1}{n(t)}dt \}=B_+(p).
$$
Similarly, we have
$$
B_-(p)=\sup_{x<m} \{ \mu((-\infty,x])\bigg[\bigg(1+\frac{(p-1)^{p/(p-2)}}{\mu(-\infty,x]}\bigg)^{(p-2)/p}-1\bigg] \int_x^m \frac{1}{n(t)} dt \}.
$$
Since $F_+'^2+F_-'^2=f'^2$ on $\rr\backslash \{m\}$, we have
$$
\aligned
&\quad \left( \int |f|^p d\mu  \right)^{2/p}-\int f^2 d\mu \\
&\le 4B_+(p)\int F_+'^2 d\nu+4B_-(p)\int F_-'^2 d\nu\\
&\le 4\max\{B_+(p),B_-(p) \}\left(\int F_+'^2 d\nu+\int F_-'^2 d\nu \right)\\
&=4\max\{B_+(p),B_-(p) \}\int f'^2 d\nu.
\endaligned
$$
Hence $C \le  4\max\{B_+(p),B_-(p) \}$.\\

{\bf Step 2.} We estimate the lower bound of $C$. At first, we suppose that $f$ is a continuous function which vanishes on $(-\infty,m]$ and is smooth on $[m,\infty)$. By approximation, $f$ satisfies (\ref{41}). Noting that in order to approach the supremum, the test function $g=-1$ on $(-\infty,m]$.   By Lemma \ref{44}, we have
$$
\aligned
&\quad \left(\int_\rr|f|^p d \mu\right)^{2/p}-\int_\rr f^2 d\mu\\
&=\sup  \left\{ \int_m^\infty f^2 g d\mu\ ;\  g\ge -1 \ \mbox{and} \ \int_m^\infty (g+1)^{p/(p-2)} d\mu \le 1 \right\}\\
&\ge \sup\left\{ \int_m^\infty f^2 g d\mu\ ;\  g\ge 0 \ \mbox{and} \ \int_m^\infty (g+1)^{p/(p-2)} d\mu \le 1 \right\}
\endaligned
$$
Since $\mu([m,\infty)) \le \frac 1 2$, we have $\int_m^\infty (g+1)^{p/(p-2)} d\mu \le 1$ for many non-negative functions $g$. By (\ref{41}), for such functions $f$ with $f(m)=0$, we have
$$
\sup\left\{ \int_m^\infty f^2 g d\mu\ ;\  g\ge 0 \ and \ \int_m^\infty (g+1)^{p/(p-2)} d\mu \le 1 \right\} \le C \int_m^\infty f'^2d \nu.
$$
By Proposition \ref{prop42}, we have
$$
C \ge \sup_{x>m} \left[\sup \left\{ \int_m^\infty 1_{[x,\infty)} g d\mu \ ;\  g \ge 0 \ and\ \int_m^\infty (g+1)^{p/(p-2)} d \mu \le 1\right\}\int_m^x \frac{1}{n(t)}dt\right].
$$
Since $\mu([m,\infty)) \le \frac 1 2 $, using Lemma \ref{lem45}, we have
$$
C \ge \sup_{x>m} \{ \mu([x,+\infty))\bigg[\bigg(1+\frac{1}{2\mu[x,+\infty)}\bigg)^{(p-2)/p}-1\bigg] \int_m^x \frac{1}{n(t)}dt \}=b_+(p).
$$
Then we suppose that  $f$ is a continuous function which vanishes on $[m, \infty)$ and is smooth on $(-\infty,m]$. Similarly, we have
$$
b_-(p)=\sup_{x<m} \{ \mu((-\infty,x])\bigg[\bigg(1+\frac{1}{2\mu(-\infty,x]}\bigg)^{(p-2)/p}-1\bigg] \int_x^m \frac{1}{n(t)}dt \}.
$$
Hence $C \ge \max \{b_+(p),b_-(p) \}$.
The proof is completed.
\nprf

\end{document}